\documentclass[12pt]{amsart}
\usepackage{amssymb,latexsym,amsmath,mathrsfs}
 \usepackage[matrix, arrow]{xy}
\xyoption{all}

\newcommand{\N}{\mathbb{N}}

\newcommand{\pf}{\noindent {\mbox{\textit{Proof}. }} }
\newcommand{\ie}{\textit{i.e.\/}\ } 

\newcommand{\cst}{C$^*$} 

\newtheorem{thm}{Theorem}[section]
\newtheorem{cor}[thm]{Corollary}

\newtheorem{prop}[thm]{Proposition}

\theoremstyle{definition}
\newtheorem{defi}[thm]{Definition}

\newtheorem{rem}[thm]{Remark}

\theoremstyle{remark}

\numberwithin{equation}{section}
\title{Extension of C$^*$-bundles}

\author{Etienne Blanchard}

\date{07/02/2008}

\begin{document}

\subjclass[2000]{Primary: 46L09; Secondary: 46L35, 46L06} 

\begin{abstract} 
We investigate which amalgamated products of continuous \cst-bundles are continuous \cst-bundles 
and we analyse the involved extension problems for continuous \cst-bundles. 
\end{abstract}

\maketitle

\section*{Introduction} 
Different (fibrewise) amalgamated products of continuous \cst-bundles have been studied 
over the last years (\cite{bla2}, \cite{KiWa}, \cite{BdWa}, \cite{bla6}), 
one of the main question being to know when 
these amalgamated products are continuous \cst-bundles. 
In order to gather these approaches in a joint framework, 
we first recall a few definitions from the theory of deformations of \cst-algebras 
and we fix several notations which will be used in the sequel. 
Then we characterise the continuity properties of different amalgamated products 
of (continuous) $C(X)$-algebras. 

\section{$C(X)$-algebras}{ 
Let $X$ be a compact Hausdorff space and 
$C(X)$ the \cst-algebra of continuous functions on $X$ 
with values in the complex field $\mathbb{C}$. 
\begin{defi} 
A $C(X)$-algebra is a \cst-algebra $A$ endowed 
with a unital $*$--homo\-morphism from $C(X)$ 
to the centre of the multiplier \cst-algebra $\mathcal{M}(A)$ of~$A$. 
\end{defi}

Given a closed subset $Y\subset X$, 
we denote by $C_0(X\setminus Y)$ the closed ideal of continuous functions on $X$ that vanish of $Y$. 
If $A$ is a $C(X)$-algebra, then the subset $C_0(X\setminus Y).A$ is a \textit{closed} ideal in $A$ 
(by Cohen factorisation Theorem) and 
we denote by $\pi^X_Y$ the quotient map $A \to A/ C_0(X\setminus Y).A\,$. 

If the closed subset $Y$ is reduced to a point $x$ and 
the element $a$ belongs to the $C(X)$-algebra $A$, 
we usually write $\pi_x$, $A_x$ and $a_x$ for $\pi^X_{\{ x\}}$, $\pi^X_{\{ x\}}(A)$ and $\pi^X_{\{ x\}}(a)$.  

Note that the function 
\begin{equation} 
x\mapsto\| \pi^X_x(a)\| =\inf\{\|\, [1-f+f(x)]a\|\,; f\in C(X)\}
\end{equation}
is always upper semi-continuous by construction. 
And the $C(X)$-algebra $A$ is said to be \textit{continuous} 
(or to be a \textit{continuous \cst-bundle over $X$}) if 
the function $x\mapsto\| \pi^X_x(a)\|$ is actually continuous for all $a$ in $A$. 

\begin{defi} 
A \textit{continuous field of states} on a unital $C(X)$-algebra $A$ is 
a unital positive $C(X)$-linear map $\varphi: A\to C(X)$. 
\end{defi}

\begin{rem} 
A (unital) separable $C(X)$-algebra $A$ is continuous if and only if (iff) 
there exists a continuous field of states $\varphi: A\to C(X)$ such that 
for all $x\in X$, the induced state $\varphi_x: a_x\in A_x\mapsto \varphi(a)(x)$ is faithful on $A_x$ 
(\cite{bla1}). 
\end{rem} 

\section{Hahn-Banach extension properties} 

Given a \cst-algebra $A$, a \cst-subalgebra $B\subset A$ 
and a state $\phi: B\to\mathbb{C}$, 
there always exists a state $\varphi$ on $A$ such that $\varphi(b)=\phi(b)$ for all $b\in B$ 
by Hahn-Banach extension theorem. 
But there is no general $C(X)$-linear version of that property. 
Indeed, consider: 
\\ \indent -- the compact space $Y:=\{ 0\}\cup \{\,\frac{1}{n}\,; n\in \mathbb{N}^*\}$, 
\\ \indent -- the unital continuous $C(Y)$-algebra $A:=C(Y)\oplus C(Y)$ and 
\\ \indent -- the $C(Y)$-subalgebra $B:=C(Y).1_A+\Bigl(C_0(Y\setminus\{0\})\oplus C_0(Y\setminus\{0\})\Bigr) \subset A$

And let $\phi: B\to C(Y)$ be the continuous field of states on $B$ fixed by the formulae 
$$
\phi(\, (b_1, b_2)\, )(\frac{1}{n}) = \left\{
\begin{array}{ll} b_1(\frac{1}{n})&\mathrm{if\; n\; is\; odd}\\
&\\ b_2(\frac{1}{n})&\mathrm{otherwise}\end{array}\right. 
\quad\mathrm{for}\quad (b_1, b_2)\in C_0(Y\setminus\{0\})\oplus C_0(Y\setminus\{0\})
$$
Then, there cannot be any continuous field of states $\varphi: A\to C(Y)$ such that 
$\varphi(b)=\phi(b)$ for all $b\in B$. 
Indeed, if $a=1\oplus 0\in A$, this would imply that 
$\varphi(a)(\frac{1}{n})=1$ if $n$ is odd and $\varphi(a)(\frac{1}{n})=0$ otherwise,  
But then, the function $y\mapsto \varphi(a)(y)$ could not be continuous at $y=0$.

\section{Tietze extension properties}
Let $(X, d)$ be a second countable compact metric space and $Y\subset X$ a non empty closed subspace. 
Given a separable continuous $C(X)$-algebra $A$, 
any continuous field of states $\phi: \pi^X_Y(A)\to C(Y)$ on the restriction $\pi^X_Y(A)$  
can always be extended to 
a continuous field of states $\varphi: A\to C(X)$ by Michael continuous selection theorem (\cite{bla1}), 
\ie such that $\varphi(a)(x)=\phi(a)(x)$ for all $a\in A$ and $x\in Y$. 

The point in this section is to study the following more general extension problem: 
Given a continuous $C(Y)$-algebra $A$, does there exist a continuous $C(X)$-algebra $D$ 
with a $C(Y)$-algebra isomorphic $\pi^X_Y(D)\cong A$? 

\medskip 
Note first that if the unital separable continuous $C(Y)$-algebra $A$ is \textit{exact}, 
then there exists a unital embedding of the $C(Y)$-algebra $A$ 
into the trivial $C(Y)$-algebra $C(Y, \otimes\mathcal{O}_2)$, 
where $\mathcal{O}_2$ is the unital Cuntz \cst-algebra generated by two isometries $s_1, s_2$ 
satisfying $1_{\mathcal{O}_2}=s_1(s_1)^*+s_2(s_2)^*$ (\cite{bla3}). 
Hence, $D:=\{\,f\in C(X, \mathcal{O}_2)\,;\, \pi^X_Y(f)\in A\,\}$ answers the question in that case. 

 \medskip 
Now, in order to study the general case, define in $X\times Y\times [0,1]$:\\
- the open subspace $U=\{ (x, y, t)\in X\times Y\times [0,1]\,;\, 0< t \}$ and \\
- the closed subspace $Z=\{ (x, y, t)\in X\times Y\times [0, 1]\,;\, 0\leq t.d(x, Y)\leq 2 d(x, Y)-d(x,y) \}$. \\ 
And let $\bar{d}$ be the metric on $Z$ given by $\bar{d}((x, y, t), (x', y', t') )=d(x, x')+d(y, y')+|t-t'|$. 

\begin{prop}\label{bw} (\cite{BdWa}) 
The coordinate map $p_1: (x, y, t)\mapsto x$ gives a structure of $C(X)$-algebra on $C(Z)$ and 
the ideal $C_0(U\cap Z)$ is  a 
continuous $C(X)$-algebra 
such that $C_0(U\cap Z)_{|Y}\cong C_0(Y\times (0, 1])$, 
i.e. the map $(x, y, t)\in U\cap Z\mapsto x\in X$ is \textit{open}. 
\end{prop} 
\pf 
Given a function $f$ in $C_0(U\cap Z)$, let us prove the continuity of the function 
$$x\in X \mapsto \|\pi_x^X(f)\|=\sup\{ |f(z)|\,;\, z\in p_1^{-1}(\{x\}) \}$$
\indent 
This map is already upper semicontinuous (u. s. c.) by construction. 
Hence, it only remains to show that for any point $x_0\in X$ and any constant $\varepsilon>0$, 
one has $\|\pi_x^X(f)\| > \|\pi_{x_0}^X(f)\|-\varepsilon$ 
for all points $x$ in a neighbourhood of $x_0$ in $X$. 

The uniform continuity of the function $f$ implies that there exists $\delta>0$ such that 
$|f(z)-f(z')|<~\varepsilon$ for all $z, z'$ in $Z$ with $\bar{d}(z, z')<\delta$. Now three cases can appear: 

1) If $x_0\in Y$ and $x\in Y$ satisfies $d(x_0, x)<\delta/2$, then $|f(x, x, t) - f(x_0, x_0, t)|<\varepsilon$ for all $t\in [0, 1]$. And so $\|\pi_x^X(f)\| > \|\pi_{x_0}^X(f)\|-\varepsilon$. 

2) If $x_0\in Y$ and $x\in X\setminus Y$ satisfies $d(x_0, x)<\delta/4$, then 
for all $y\in Y$, the relation $d(x, y)\leq 2 d(x, Y)$ implies that 
$d(y, x_0)\leq d(y, x)+d(x, x_0)\leq 2 d(x, Y)+ d(x, x_0)\leq \frac{3}{4}\delta$ and so 
$\bigl| f(x, y, t) - f(x_0, x_0, t)\bigr|<\varepsilon$ for all $t\in [0, 2 -\frac{d(x, y)}{d(x, Y)} ]$. 
Whence the inequality $\|\pi_x^X(f)\| > \|\pi_{x_0}^X(f)\|-\varepsilon$. 

3) If $x_0\not\in Y$ and the triple $(x_0, y_0, t_0)\in U\cap Z$ satisfies 
$| f(x_0, y_0, t_0) | = \|\pi_{x_0}^X(f)\| \not= 0$, then $d(x_0, y_0)< 2 d(x_0, Y)$. 
Thus, there exists by continuity a constant $\alpha(x_0)\in ] 0, \delta/2[$ such that 
all $x\in X$ in the ball of radius $\alpha(x_0)$ around $x_0$ satisfy: 
\begin{center}
a) $d(x, Y)>0\,$, \quad 
b) $d(x, y_0) <2 d(x, Y)\,$, \quad 
c) $t_0 < 2-\frac{d(x, y_0)}{d(x, Y)} + \delta/2\,.$
\end{center}
And so $\|\pi_x^X(f)\| \geq \left|f(x, y_0, \inf\{ t_0, 2-\frac{d(x, y_0)}{d(x, Y)}\})\right| > 
\|\pi_{x_0}^X(f)\|-\varepsilon$. 
\qed
\begin{rem} S. Wassermann pointed out that if $Y=\{ 0, 1\}\subset X=[0,1]$, 
then $Z=\{(x, 0, t)\in [0, 1]\times \{0\}\times [0,1]\,;\, t\leq\frac{2-3x}{1-x}\} \cup
\{(x, 1, t)\in [0, 1]\times \{1\}\times [0,1]\,;\, t\leq\frac{3x-1}{x}\}$. 
Hence, this $C(X)$-algebra $C(Z)$ is not continuous at $x=\frac{1}{3}$ and $x=\frac{2}{3}$. 
\end{rem}

Above Proposition \ref{bw} implies the following. 
\begin{cor}\label{extension} 
Let $A$ be a continuous $C(Y)$-algebra.
\begin{enumerate}
\item[a)] $B:= C(X)\otimes A\otimes C([0,1])$ is a continuous $C(X\times Y\times [0, 1])$-algebra. 
\item[b)] $D:=\left[ C_0(U). B\right]_{|Z}=C_0(U).B \big/ C_0(U\backslash U\cap Z).B$ is a continuous $C(X)$-algebra. 
\item[c)] There is an isomorphism of $C(Y)$-algebras $D_{|Y}\cong A\otimes C_0( (0,1] )$.  
\end{enumerate}
\end{cor}
\pf 
b) Let $b\in D$. Then for all $x\in X$, we have  
$$\|\pi_x^X(b)\|=\| b+C_0(X\setminus\{x\})D\|=\sup\{\| \pi_z^Z(b) \|\,;\, z\in  p_1^{-1}(\{ x\})\}\,,$$
whence the continuity of the map $x\mapsto\|\pi_x^X(b)\|$ by a) and Proposition~\ref{bw}. \qed

\section{Amalgamated tensor products of continuous $C(X)$-algebras}
Given a fixed compact Hausdorff space $X$, we study in this section the continuity properties of 
the different tensor products amalgamated over $C(X)$ of 
two given continuous $C(X)$-algebras $A$ and $B$. 

Let $A\odot B$ denote the algebraic tensor product (over $\mathbb{C}$) of $A$ and $B$, 
let $\mathcal{I}_X(A, B)$ be the ideal in $A\odot B$ generated by the differences 
$a f\otimes b - a\otimes f b$ ($a\in A$, $b\in B$, $f\in C(X)\,$) 
and let $A{\mathop{\odot}\limits_{C(X)} B}$ denote the quotient of $A\odot B$ by $\mathcal{I}_X(A, B)$. 
\\ \indent 
If $C_\Delta(X\times X)\subset C(X\times X)$ is the ideal of continuous function of $X\times X$ 
which are zero on the diagonal and 
$A\mathop{\otimes}\limits^m B$ (resp. $A\mathop{\otimes}\limits^M B$) 
is the minimal (resp. maximal) tensor product over $\mathbb{C}$ 
of the two \textit{continuous} $C(X)$-algebras $A$ and $B$, 
then the quotient $A\mathop{\otimes}\limits^m_{C(X)} B:=
A\mathop{\otimes}\limits^m B\Big{\slash}C_\Delta(X\times X) A\mathop{\otimes}\limits^m B$ 
(resp. $A\mathop{\otimes}\limits^M_{C(X)} B:=
A\mathop{\otimes}\limits^M B\Big{\slash}C_\Delta(X\times X) A\mathop{\otimes}\limits^M B$) 
is the minimal (resp. maximal) completion of the algebraic amalgamated tensor product 
$A\mathop{\odot}\limits_{C(X)} B$. 
Further, the $*$-algebra $A\mathop{\odot}\limits_{C(X)} B$ embeds 
in the $C(X)$-algebra $A\mathop{\otimes}\limits^m_{C(X)} B$ (\cite{bla2}) and we have 
\begin{equation} 
\forall\,x\in X\,,\hspace{20pt}  
(A\mathop{\otimes}\limits^m_{C(X)} B)_x\cong A_x\mathop{\otimes}\limits^m B_x
\hspace{20pt} \mathrm{and}\hspace{20pt} 
(A\mathop{\otimes}\limits^M_{C(X)} B)_x\cong A_x\mathop{\otimes}\limits^M B_x\,. 
\end{equation} 

Let us also recall a characterisation of exactness given by Kirchberg and Wassermann. 
\begin{prop}\label{exactness-crit}(\cite[Theorem 4.5]{KiWa}) 
Let $Y=\mathbb{N}\cup\{\infty\}$ be the one point compactification of $\mathbb{N}$ 
and let $D$ be a \cst-algebra. 
Then the following assertions are equivalent. 
\begin{enumerate} 
\item[i)] The \cst-algebra $A$ is exact. 
\item[ii)] For all continuous $C(Y)$-algebra $B$, the minimal tensor product $A\mathop{\otimes}\limits^m B$ is a continuous $C(Y)$-algebra with fibres $A\mathop\otimes\limits^m B_y$ ($y\in Y$). 
\end{enumerate} 
\end{prop}

Then, the following holds. 
\begin{prop}(\cite{BdWa}, \cite{bla6}) 
Let $X$ be a second countable compact Hausdorff space and 
$A$ a separable unital continuous $C(X)$-algebra. 
\\ \indent 
If the topological space $X$ is \textit{perfect} (\ie without isolated point), 
then the following assertions $\alpha_e$ and $\gamma_e$ 
(resp. $\alpha_n$ and $\gamma_n$) are equivalent. 

\noindent 
$\alpha_e$) The \cst-algebra $A$ is exact. 

\noindent 
$\gamma_e$) For all continuous $C(X)$-algebra $D$, 
the amalgamated tensor product $\mathrel{A\mathop{\otimes}\limits^m_{C(X)} D}$ is a 
continuous $C(X)$-algebra with fibres $\mathrel{A_x\mathop{\otimes}\limits^m D_x}$ ($x\in X$). 

\noindent 
$\alpha_n$) The \cst-algebra $A$ is nuclear. 

\noindent 
$\gamma_n$) For all continuous $C(X)$-algebra $D$, 
the amalgamated tensor product $\mathrel{A\mathop{\otimes}\limits^M_{C(X)} D}$ is a 
continuous $C(X)$-algebra with fibres $\mathrel{A_x \mathop{\otimes}\limits^M D_x}$ ($x\in X$). 
\end{prop}
\pf 
$\alpha_e)\Rightarrow \gamma_e)$ If the \cst-algebra $A$ is exact, 
then $A\mathop{\otimes}\limits^m D$  is a continuous $C(X\times X)$-algebra 
with fibres $A_x\mathop{\otimes}\limits^m D_{x'}$ ($x, x'\in X$) (\cite{KiWa}). 
Hence, its restriction to the diagonal is as desired. 

\medskip 
\noindent $\gamma_e)\Rightarrow \alpha_e)$ 
Suppose conversely that the $C(X)$-algebra $A$ satisfies $\gamma_e$). 

\noindent a) All the fibres $A_x$ are exact ($x\in X$). 
Indeed, given a point $x$ in $X$, take a sequence of points $x_n$ in $X$ converging to $x$ 
such that there is a topological isomorphism 
$Y:=\{x_n; n\in\N\}\cup\{ x\}\cong \mathbb{N}\cup\{\infty\}$. 
Then, for any separable continuous $C(Y)$-algebra $D$, there is  a continuous $C(X)$-algebra 
$\mathcal{D}$ such that $\mathcal{D}_{|Y}= D\otimes C_0( (0, 1])$ (Corollary~\ref{extension}). 
Now, the continuity of the $C(X)$-algebra $\mathcal{D}\mathop{\otimes}\limits^m_{C(X)} A$ 
given by $\gamma_e$ 
implies that of its restriction $\left(\mathcal{D}\mathop{\otimes}\limits^m_{C(X)} A\right)_{|Y}\cong
\left( C_0( (0,1])\otimes D\right)\mathop{\otimes}\limits^m_{C(Y)} A_{|Y}$, 
whence that of the $C(Y)$-algebra $D \mathop{\otimes}\limits^m_{C(Y)} A_{|Y}$ 
since there is an isometric $C(Y)$-linear embedding $D\hookrightarrow \mathcal{D}_{|Y}$.  
And this implies the exactness of the \cst-algebra $A_x$ is exact by Proposition~\ref{exactness-crit}. 

\noindent \textit{b)} If $B$ is a \cst-algebra and 
$\mathcal{B} $ is the constant $C(X)$-algebra $C(X; B)$, 
then for all $x\in X$, we have the exact sequence \\
\centerline{$0\to C_x(X) A\mathop\otimes\limits^m B\to 
(A\mathop{\otimes}\limits^m_{C(X)} \mathcal{B} )_x=A\mathop\otimes\limits^m B\to 
A_x\mathop\otimes\limits^m B\to 0$.} 

\noindent \textit{c)} 
If $D$ is a $C(X)$-algebra, then for all point $x\in X$, 
we have the sequence of epimorphisms 
$(A\mathop{\otimes}\limits_{C(X)}^m D)_x\twoheadrightarrow 
(A_x\mathop{\otimes}\limits^m D)_x \twoheadrightarrow 
A_x\mathop{\otimes}\limits^m D_x$

\noindent \textit{d)} Now, let $B$ be a \cst-algebra, $K\triangleleft B$ a closed two sided ideal in $B$ 
and take an element $d\in\ker\{ A\mathop\otimes\limits^m B\to A\mathop\otimes\limits^m B/K\}$. 
Then for all $x\in X$, we have 
$$\begin{array}{clr}
d_x\!\!\!&\in\ker\{ (A\mathop\otimes\limits^m B)_x\to (A\mathop\otimes\limits^m B/K)_x\}&\\
&=\ker\{ A_x\mathop\otimes\limits^m B \to A_x\mathop\otimes\limits^mB/K\} 
&\quad\mathrm{by}\, \textit{b)}\\ 
&=A_x\mathop\otimes\limits^m K &\quad \mathrm{by}\, \textit{a)}\\
&=(A\mathop\otimes\limits^m K)_x 
&\quad \mathrm{by}\, \textit{c)}
\end{array}$$
Thus, $d\in A\mathop\otimes\limits^m K$. 
And so, the \cst-algebra $A$ is exact. 

\medskip 
The proof of $\alpha_n)\Rightarrow \gamma_n)$ is similar to that of $\alpha_e)\Rightarrow \gamma_e)$. 
On the other hand, if a \cst-algebra $A$ satisfies $\gamma_n$), 
then all the fibres $A_x$ ($x\in X$) are nuclear by \cite[Theorem~3.2]{KiWa} 
and so the \cst-algebra $A$ itself is nuclear (see e.g. \cite[Proposition 3.23]{bla1}). 
\qed

\begin{rem} 
These characterisations do not hold anymore if the compact space $X$ is not perfect. 
Indeed, if the space $X$ is reduced to a point, then both the amalgamated tensor products 
$A\mathop{\otimes}\limits^m_{C(X)} D$ and $A\mathop{\otimes}\limits^M_{C(X)} D$ 
are constant, hence continuous. 
\end{rem}

\section{Amalgamated free products of continuous $C(X)$-algebras}
We now describe the continuity properties of different free products amalgamated over $C(X)$ 
of two given unital continuous $C(X)$-algebras $A$ and $B$.

\begin{prop}(\cite[Corollary 4.8]{bla6}) 
Let $X$ be a second countable perfect compact Hausdorff space and 
$A$ a separable unital continuous $C(X)$-algebra. 
Then the following assertions are equivalent. 

\noindent 
$\alpha_e$) The \cst-algebra $A$ is exact. 

\noindent 
$\gamma_e$) For all separable unital continuous $C(X)$-algebra $D$ and 
all continuous fields of states $\phi: A\to C(X)$, $\psi: D\to C(X)$, 
the $C(X)$-algebra $(A, \phi)\mathop{\ast}\limits^r_{C(X)} (D, \psi)$ is continuous. 
\end{prop}

\begin{rem}
The is no similar result for full amalgamated free product. 
Indeed, the full amalgamated free product $A\mathop{\ast}\limits^f_{C(X)} D$ 
of two unital continuous $C(X)$-algebras $A$ and $D$ is always a continuous $C(X)$-algebra 
with fibres $A_x\mathop{\ast}\limits^f D_x$ ($x\in X$) 
(\cite[Theorem~3.7]{bla6}). 
\end{rem}

\noindent
\email{Etienne.Blanchard@math.jussieu.fr}\\
\address{IMJ,
175, rue du Chevaleret, F--75013 Paris}

\end{document}